\magnification=\magstep1
\input amstex

\input prepictex
\input pictex
\input postpictex
\documentstyle{amsppt}
\TagsOnRight
\hsize=5in
\vsize=7.8in
\def\buildrel#1\over#2{\mathrel{\mathop{\null#2}\limits^{#1}}}
\def\buildrul#1\under#2{\mathrel{\mathop{\null#2}\limits_{#1}}}
\topmatter

\title Witten deformation and \\
polynomial differential forms\endtitle
\author Michael Farber and Eugenii Shustin\endauthor
\address
School of Mathematical Sciences,
Tel-Aviv University,
Ramat-Aviv 69978, Israel
\endaddress
\email farber\@math.tau.ac.il, shustin\@math.tau.ac.il
\endemail
\thanks The first author was supported by a grant from the
US - Israel Binational Science Foundation. The second author was
partially supported by Grant
No. 6836-1-9 of the Ministry of Science and Technology, Israel.
\endthanks
\abstract{As is well-known, the Witten deformation $d_h$ 
of the De Rham complex computes the De Rham cohomology.
In this paper we study the Witten deformation on a noncompact manifold 
$M$ and restrict it to differential
forms which behave {\it polynomially near infinity}.
Such polynomial differential forms naturally appear on manifolds with {\it a cylindrical structure}. 
We prove that the cohomology of the Witten deformation $d_h$ acting on the complex of the
polynomially growing forms (depends on $h$ and) 
can be computed as the relative cohomology of the pair $(M,F)$ 
where $F$ 
is a {\it negative remote fiber of } $h$. We show that the assumptions of our main theorem are satisfied in a number of interesting special cases, including 
generic real polynomials
on ${\bold R}^n$.}
\endabstract
\endtopmatter
\def\figfirst{%
$$
\font\thinlinefont=cmr5
\beginpicture
\setcoordinatesystem units < 0.900cm, 0.900cm>
\linethickness=0pt
\putrule from -7.000 0.000  to  14.000 0.000
\putrule from -7.000 0.000  to  -7.000 6.000
\setshadesymbol ({\thinlinefont .})
\setlinear
%
%
\linethickness=1pt
\setplotsymbol ({\thinlinefont .})
\circulararc 90 degrees from  0.000 4.000 center at  -3.000 1.000
\circulararc 90 degrees from  -6.000 4.000 center at  -5.000 3.000
\circulararc 90 degrees from  -6.000 2.000 center at  -3.000 5.000
\circulararc 53.130 degrees from  -2.000 3.000 center at  -3.000 1.000
\circulararc 90 degrees from  -4.581 3.419 center at  -3.000 5.000
\circulararc 45 degrees from  0.000 4.000 center at  1.000 5.000
\circulararc 45 degrees from  1.000 2.414 center at  1.000 1.000
\circulararc 60 degrees from  3.500 3.586 center at  3.500 3.500
\circulararc 60 degrees from  3.425 3.543 center at  4.366 3.000
\circulararc 60 degrees from  3.425 2.457 center at  3.500 2.500
\linethickness=0.5pt
\putrule from  1.000 2.414  to  6.000 2.414
\putrule from  1.000 3.586  to  6.000 3.586
\putrule from  1.500 3.386  to  3.000 3.386
\putrule from  1.500 2.614  to  3.000 2.614
\putrule from  1.250 3.000  to  2.750 3.000
\putrule from  4.500 3.386  to  5.500 3.386
\putrule from  4.500 2.614  to  5.500 2.614
\putrule from  4.250 3.000  to  5.250 3.000
\setdashes
\linethickness=1pt
\circulararc 60 degrees from  3.500 2.414 center at  3.500 2.500
\circulararc 60 degrees from  3.575 2.457 center at  2.634 3.000
\circulararc 60 degrees from  3.575 3.543 center at  3.500 3.500
\circulararc 60 degrees from  6.000 2.414 center at  6.000 2.500
\circulararc 60 degrees from  6.075 2.457 center at  5.134 3.000
\circulararc 60 degrees from  6.075 3.543 center at  6.000 3.500
\circulararc 60 degrees from  6.000 3.586 center at  6.000 3.500
\circulararc 60 degrees from  5.925 3.543 center at  6.866 3.000
\circulararc 60 degrees from  5.925 2.457 center at  6.000 2.500

\endpicture
$$%
}

\def\figthird{%
$$
\font\thinlinefont=cmr5
\beginpicture
\setcoordinatesystem units < 0.900cm, 0.900cm>
\linethickness=0pt
\putrule from -7.000 0.000  to  14.000 0.000
\putrule from -7.000 0.000  to  -7.000 6.000
\setshadesymbol ({\thinlinefont .})
\setlinear
%
%
\linethickness=1pt
\setplotsymbol ({\thinlinefont .})
\circulararc 90 degrees from  0.000 4.000 center at  -3.000 1.000
\circulararc 90 degrees from  -6.000 4.000 center at  -5.000 3.000
\circulararc 90 degrees from  -6.000 2.000 center at  -3.000 5.000
\circulararc 53.130 degrees from  -2.000 3.000 center at  -3.000 1.000
\circulararc 90 degrees from  -4.581 3.419 center at  -3.000 5.000
\circulararc 45 degrees from  0.000 4.000 center at  1.000 5.000
\circulararc 45 degrees from  1.000 2.414 center at  1.000 1.000
\circulararc 60 degrees from  2.500 3.586 center at  2.500 3.500
\circulararc 60 degrees from  2.425 3.543 center at  3.366 3.000
\circulararc 60 degrees from  2.425 2.457 center at  2.500 2.500
\plot  3.450 5.200  3.500 5.500  3.550 5.200 /
\plot  1.000 0.750  3.500 2.000 /
\plot  1.500 4.000  3.500 5.000 /
\plot  1.243 0.926  1.000 0.750  1.286 0.838 /
\plot  5.200 2.050  5.500 2.000  5.200 1.950 /
\linethickness=0.5pt
\putrule from  1.000 2.414  to  1.500 2.414
\putrule from  6.000 2.414  to  2.500 2.414
\putrule from  1.000 3.586  to  1.500 3.586
\putrule from  6.000 3.586  to  2.500 3.586
\putrule from  3.500 2.000  to  3.500 2.414
\putrule from  3.500 5.500  to  3.500 3.586
\putrule from  1.500 1.000  to  1.500 4.000
\putrule from  3.500 2.000  to  5.500 2.000
\setdashes
\linethickness=1pt
\circulararc 60 degrees from  2.500 2.414 center at  2.500 2.500
\circulararc 60 degrees from  2.575 2.457 center at  1.634 3.000
\circulararc 60 degrees from  2.575 3.543 center at  2.500 3.500
\circulararc 60 degrees from  6.000 2.414 center at  6.000 2.500
\circulararc 60 degrees from  6.075 2.457 center at  5.134 3.000
\circulararc 60 degrees from  6.075 3.543 center at  6.000 3.500
\circulararc 60 degrees from  6.000 3.586 center at  6.000 3.500
\circulararc 60 degrees from  5.925 3.543 center at  6.866 3.000
\circulararc 60 degrees from  5.925 2.457 center at  6.000 2.500
\put {$x_1$} at 1.400 0.650
\put {$x_{n-1}$} at 4.050 5.450
\put {$x_n=t$} at 5.400 1.750

\endpicture
$$%
}
\def\figsecond{%
$$
\font\thinlinefont=cmr5
\beginpicture
\setcoordinatesystem units < 0.900cm, 0.900cm>
\linethickness=0pt
\putrule from -7.000 0.000  to  14.000 0.000
\putrule from -7.000 0.000  to  -7.000 6.000
\setshadesymbol ({\thinlinefont .})
\setlinear
\linethickness=1pt
\putrule from  -6.000 1.242  to  6.000 1.242
\putrule from  -6.000 4.758  to  6.000 4.758
\circulararc 60 degrees from  1.500 4.758 center at  1.500 4.500
\circulararc 60 degrees from  1.275 4.629 center at  4.098 3.000
\circulararc 60 degrees from  1.275 1.371 center at  1.500 1.500
\circulararc 60 degrees from  -1.500 4.758 center at  -1.500 4.500
\circulararc 60 degrees from  -1.725 4.629 center at  1.098 3.000
\circulararc 60 degrees from  -1.725 1.371 center at  -1.500 1.500
\ellipticalarc axes ratio 8:1 -175.5 degrees from  -6.432 4.200 center at  -6.700 3.700
\ellipticalarc axes ratio 7:1 -169 degrees from  5.568 1.800 center at  5.800 2.300
\linethickness=0.6pt
\putrule from  -6.345 1.600  to  5.655 1.600
\putrule from  -6.345 4.400  to  5.655 4.400
\putrule from  -6.661 3.000  to  5.339 3.000
\linethickness=0.3pt
\plot 0.000 5.200  0.000 3.050 /
\plot -0.050 3.550  0.000 3.050  0.050 3.550 /
\plot -0.400 5.200  -0.700 4.450 /
\plot 0.400 5.200  0.700 1.650 /
\plot -0.670 4.660  -0.700 4.450  -0.580 4.640 /
\plot 0.710 2.055  0.700 1.650  0.625 2.050 /

\setdashes
\linethickness=0.5pt
\circulararc 60 degrees from  1.500 1.242 center at  1.500 1.500
\circulararc 60 degrees from  1.725 1.371 center at  -1.098 3.000
\circulararc 60 degrees from  1.725 4.629 center at  1.500 4.500
\circulararc 60 degrees from  -1.500 1.242 center at  -1.500 1.500
\circulararc 60 degrees from  -1.275 1.371 center at  -4.098 3.000
\circulararc 60 degrees from  -1.275 4.629 center at  -1.500 4.500
\circulararc 60 degrees from  6.000 1.242 center at  6.000 1.500
\circulararc 60 degrees from  6.225 1.371 center at  3.402 3.000
\circulararc 60 degrees from  6.225 4.629 center at  6.000 4.500
\circulararc 60 degrees from  -6.000 1.242 center at  -6.000 1.500
\circulararc 60 degrees from  -5.775 1.371 center at  -8.598 3.000
\circulararc 60 degrees from  -5.775 4.629 center at  -6.000 4.500
\circulararc 60 degrees from  6.000 4.758 center at  6.000 4.500
\circulararc 60 degrees from  5.775 4.629 center at  8.598 3.000
\circulararc 60 degrees from  5.775 1.371 center at  6.000 1.500
\circulararc 60 degrees from  -6.000 4.758 center at  -6.000 4.500
\circulararc 60 degrees from  -6.225 4.629 center at  -3.402 3.000
\circulararc 60 degrees from  -6.225 1.371 center at  -6.000 1.500

\put {$t=1$} at 1.500 0.700
\put {$t=-1$} at -1.500 0.700
\put {$tg(x)\le-c$} at -4.100 3.700
\put {$tg(x)\le-c$} at 4.000 2.300
\put {$g(x)=0$} at 0.000 5.500
\put {$M=N\times\R$} at 0.000 0.000
\endpicture
$$%
}

\def\figfourth{%
$$
\font\thinlinefont=cmr5
\beginpicture
\setcoordinatesystem units < 0.900cm, 0.900cm>
\setlinear
\linethickness=1pt
\putrule from -5.000 3.000  to  5.000 3.000
\putrule from 0.000 0.000  to  0.000 6.000
\plot  4.750 3.050  5.000 3.000  4.750 2.950 /
\plot  -0.050 5.750  0.000 6.000  0.050 5.750 /
\setquadratic \plot  4.500 5.419  2.000 3.000  4.500 0.581 /
\plot  -4.500 5.419  -2.000 3.000  -4.500 0.581 /
\setlinear
\linethickness=0.3pt
\plot  -1.020 1.400  -0.700 2.800 /
\plot  -0.720 2.520  -0.700 2.800  -0.800 2.530 /
\plot 2.200 3.200  3.900 4.900 /
\plot 2.700 3.200  4.500 5.000 /
\plot 3.200 3.200  4.600 4.600 /
\plot 3.700 3.200  4.700 4.200 /
\plot 2.400 2.800  2.200 2.600 /
\plot 2.900 2.800  2.350 2.250 /
\plot 3.400 2.800  2.600 2.000 /
\plot 3.900 2.800  2.850 1.750 /
\plot 4.400 2.800  3.150 1.550 /
\plot 4.900 2.800  3.450 1.350 /
\plot 4.850 2.250  3.750 1.150 /
\plot 4.800 1.700  4.050 0.950 /
\plot 4.750 1.150  4.350 0.750 /

\plot -2.200 2.800  -3.900 1.100 /
\plot -2.700 2.800  -4.500 1.000 /
\plot -3.200 2.800  -4.600 1.400 /
\plot -3.700 2.800  -4.700 1.800 /
\plot -4.200 2.800  -4.800 2.200 /
\plot -4.700 2.800  -4.900 2.600 /

\plot -2.400 3.200  -2.200 3.400 /
\plot -2.900 3.200  -2.350 3.750 /
\plot -3.400 3.200  -2.600 4.000 /
\plot -3.900 3.200  -2.850 4.250 /
\plot -4.400 3.200  -3.150 4.450 /
\plot -4.900 3.200  -3.450 4.650 /
\plot -4.850 3.750  -3.750 4.850 /
\plot -4.800 4.300  -4.050 5.050 /
\plot -4.750 4.850  -4.350 5.250 /
\put {$M_+$} at  0.500 5.750
\put {$M_-$} at  4.750 3.250
\put {$D(M_-)$} at  -1.000 1.000
\linethickness=2.5pt
\putrule from  -2.500 3.000  to  2.500 3.000
\endpicture
$$%
}

\define\C{{\bold C}}
\define\R{{\bold R}}

\define\T{{\Cal T}}

\define\M{{\Cal M}}
\redefine\P{{\Cal P}}

\define\sign{\operatorname{sign}}


\define\id{\operatorname{id}}

\define\supp{\operatorname{supp}}


\def\<{\langle}
\def\>{\rangle}

\define\pd#1#2{\dfrac{\partial#1}{\partial#2}}

\define\cc{\Cal C}

\define\Cone{\operatorname{Cone}}
\redefine\gg{\frak g}
\define\p{\Cal P}
\define\po{\p\Omega}

\def\part{\partial}

\NoBlackBoxes
\documentstyle{amsppt}

\nopagenumbers

\heading{\bf \S 0. Introduction}\endheading

E. Witten, in his paper on the Morse theory \cite{W},
invented the following deformed differential
$$d_{th} = d + tdh\wedge\cdot,\quad t\in \R,$$
acting on the De Rham complex of smooth differential forms $\Omega^\bullet(M)$ on a compact
manifold $M$; here $h: M\to \R$ is a Morse function. The main point of Witten's approach was the
observation that, since $d_{th}(\omega) = e^{-th}d(e^{th}\omega)$,
the deformed differential is gauge equivalent to the exterior derivative $d$ and so produces
the same cohomology. However for large $t$ the structure of the small eigenfunctions of the corresponding Laplacian depends on $h$ and can be described in terms of the critical points of $h$.

The Theorem of A. Dimca and M. Saito \cite{DS} claims that if $h:\C^n\to \C$ is a complex polynomial
and if one considers the deformed differential
$d_h = d + dh\wedge\cdot$ acting on the space of polynomial differential forms on $\C^n$, then
the cohomology in dimension $k$ equals the reduced cohomology $\tilde H^{k-1}(F)$, where
$F=h^{-1}(z)$, $z\in \C$ is a generic fiber of $h$. In this situation the gauge transformation
$\omega\mapsto e^h\omega$ is not relevant since it does not preserve the space of polynomial
forms.

Alexander and Maxim Braverman in \cite{BB} suggested a real version of the Dimca-Saito theorem.
They considered the situation, when $h:\R^n\to \R$ is a real polynomial and the Witten deformed
differential $d_h$ acts on the space of tempered currents on $\R^n$. Their result states that
under certain conditions on $h$ (for example, if $h$ is homogeneous) the cohomology of the
De Rham complex of tempered currents with respect to $d_h$ equals $\tilde H^{k-1}(F)$,
where $F$ is a "negative remote fiber of $h$"; in other words, $F=h^{-1}(-c)$ where $c>0$ is large
enough. Note that the tempered currents are those,  which behave {\it polynomially near infinity.}
In \cite{BB} A. and M. Braverman conjectured that their theorem must hold true for any real polynomial.

The purpose of the present paper is to compute the cohomology of the Witten deformation
$d_h$ for noncompact manifolds $M$, having {\it cylindrical structure near infinity}. 
We show that for manifolds with cylindrical structure one naturally defines {\it the complex of polynomially growing forms} $\po^\bullet(M)$, which is a subcomplex of the
De Rham complex. Under some natural assumptions the Witten deformed differential $d_h$ 
acts on the complex of polynomially growing forms. Our main result computes
the cohomology of $(\po^\bullet(M), d_h)$, as the relative cohomology of $(M,F)$, where
$F$ is a remote negative fiber of $h$.

Formally, the definition of cylindrical structure involves two components: the first is a smooth
function $t:M\to \R$ (which intuitively measures distances to the compact part of $M$);
the second is a certain Lie algebra of vector fields $\gg$ on $M$, which we call {\it constant
vector fields}. Intuitively, the constant vector fields
determine the {\it differential scale near infinity}. We define polynomially growing functions on
$M$ as those $f:M\to \R$ such that all their higher derivatives $X_1X_2\dots X_k(f)$
admit polynomial estimates in $t$.

Our first main theorem (Theorem 1) gives conditions under which cohomology of the complex
of polynomially growing forms coincides with the relative cohomology of the remote negative 
fiber. Here we use the notion of {\it vector field development} which allows one to control the estimates of the higher
derivatives of the diffeomorphisms generated by a vector field.
We apply our general Theorem 1 to some special situations which arise in applications to Morse theory (Theorems 2 and 3).

$\R^n$ has many different cylindrical structures. Here it is natural to consider the function
$r=r(x)=|x|, \quad 
x\in \R^n$ as the function $t$, and there are different inequivalent choices
for $\gg$. One such choice (called the {\it standard cylindrical structure on $\R^n$})
 consists of choosing  $\gg$ as the set of vector fields commuting with the Euler (radial) vector field
$\partial_r=r^{-1}\sum x_i\partial_i$.
We prove that if $h:\R^n\to \R$ is a real generic polynomial then the cohomology of the Witten deformation $d_h$ acting on the complex of
polynomially growing forms (with respect to the standard cylindrical structure) coincides with
the cohomology of the remote negative fiber. More precisely, we prove the above statement
in several different cases: 

A. if $h$ is quasi-homogeneous;

B. the Newton polyhedron of $h$ lies under a hyperplane
$\sigma$, and the polynomial $h^{\sigma}$
has no critical points on $\R^n$ except the origin;

C. $h$ is bounded from below.

Note that the case B includes a generic polynomial.

We also prove here an interesting theorem (Theorem 5), stating that any real polynomial near infinity 
is conjugate to its principal
leading part (assuming that the latter is non-degenerate) and the intertwining
diffeomorphism may be chosen so that all its partial derivatives (of all orders) are bounded.

The authors are grateful to M. Braverman for useful discussions.

\heading{\bf \S 1. Statement of the main results}\endheading

\subheading{1.1. Cylindrical structure} Let $M$ be a noncompact manifold
with a cylindrical end (see Fig. 1).
This means that $M$ is represented as a union
$$M = M_0\cup \Cal T,\quad M_0\cap \Cal T = \partial M_0=\partial \Cal T,\tag1-1$$
where $M_0$ is a compact manifold with boundary $\partial M_0$, and $\Cal T$
is represented as the product
$$\Cal T =  \partial M_0\times [1,\infty)\tag1-2$$
with $\partial M_0 \times 1$ identified with $\partial M_0$.
$\T$ will be called {\it "the tube part of $M$"}. We have the natural coordinates $(x,t)$ on $\T$,
where $x\in \partial M_0$ and $t\in [1,\infty)$.

Representations (1-1) and (1-2) determine {\it the cylindrical structure in the neighborhood of the end}.

Consider simple examples of manifolds with cylindrical ends.

1. $\R^n$. Here the compact part $M_0$ consists of the unit ball $B\subset \R^n$ and the coordinates
on the tube part $\T = \R^n- B$ are given by $r \mapsto (|r|^{-1}\cdot r, |r|)$.

2. Let $E$ be the total space of a smooth vector bundle over a compact manifold $N$. We can
define the cylindrical structure on $E$ by taking for the compact part the fibration on unit balls
over $N$, and by identifying the complement with the product of the sphere bundle with half
line. This is a parametrized version of the previous example.

3. Any compact manifold with boundary $M$
can be viewed as a manifold with a cylindrical end as follows.
A neighborhood of the boundary can be represented as $\partial M \times [0,1)$, where
$\partial M = \partial M\times {0}$. Thus the open manifold $M - \partial M$ has the tube part
$\partial M \times (0,1)$.

4. An interesting class of examples of cylindrical structures can be constructed as the complements of divisors.
We will mention here only the case of simple divisors formed by smooth submanifolds; more interesting examples yielding divisors with normal crossings, will be considered elsewhere.
Let $\M$ be a smooth closed manifold and let $\Sigma\subset \M$ be its smooth closed submanifold.
We will consider the complement
$$M=\M-\Sigma$$
with the following cylindrical structure. Choose
a Riemannian metric on $\M$ and let
$$t=(r_\Sigma)^{-2}: M\to \R,$$
where $r_\Sigma(x)$ denotes the
distance from $x\in \M$ to $\Sigma$. Note that $t$ is smooth near $\Sigma$, i.e. on a 
complement of a compact subset in $M$. To construct a
vector field $\partial_t$, we construct the following vector field 
$Y$ on a neighborhood of
infinity in $M$. Given a point $x\in M$, which is close enough to $\Sigma$ in $\M$, we may join
it to $\Sigma$ by a unique shortest geodesic of length $r_\Sigma(x)$; we 
will let $Y_x$ to be
the unit vector tangent to this geodesic at $x$. This defines a smooth 
vector field $Y$
on a complement of a compact subset of $M$ so that $Y(r_\Sigma)=-1$. Now we set
$$\partial_t =1/2\cdot t^{3/2}Y.$$
The pair $(t, \partial_t)$ provides a cylindrical structure on $M = \M-\Sigma$.
\topinsert
\figfirst
\botcaption{Figure 1}\endcaption
\endinsert

Our aim now is to give a precise definition of a cylindrical structure.
Suppose that $M$ is a non-compact manifold, $t:M\to \R$
is a function, which is smooth on the complement of a compact subset of $M$,
and $X=\partial_t$ is a smooth vector field on $M$ with the following properties:
\roster
\item $t:M\to \R$ is a proper map with image $t(M)$ 
lying in an interval of the form $[a,\infty)\subset \R$.
\item
$\partial_t(t)= 1$ holds on a COCS in $M$, cf. below.
\endroster
In this paper we will use the following convention:
$$\text{COCS= complement of a compact subset.}$$

The pair $(t,\partial_t)$ as above 
determines a representation (1-1). Namely, we will consider the fibers
$t^{-1}(c)$ for large $c\in\R$ and the flow on $M$,
determined by $\partial_t$, will represent the complement of a compact subset of $M$ as the
product $t^{-1}(c)\times \R_+$.

\proclaim{Definition} A choice of cylindrical structure on $M$ consists of specifying of a pair
$(t, \partial_t)$ as above.\endproclaim

\subheading{1.2. Constant vector fields}
We will see that in general a choice of cylindrical structure leads to a few natural algebras of vector fields on $M$ (which we call {\it algebras of constant, bounded and polynomial vector fields}) 
and to a subcomplex of polynomially growing differential forms on $M$.

Assume that a cylindrical structure is given by a pair $(t,\partial_t)$ as above.
A smooth vector field $Y$ on $M$ will be called {\it constant} if it commutes
with $\partial_t$ on a COCS, i.e. if the support of $[Y,\partial_t]$ is compact.
We will denote the set of all constant vector fields  by $\gg$. It is clearly a Lie algebra.

Denote by $\cc=\cc(M)$ the set of all functions $f:M\to \R$, such that the support of $\partial_t(f):M\to \R$ is compact.
$\cc$ is a ring. Functions $f\in\cc$ are eventually constant along the trajectories of $\partial_t$. It is clear that
multiplying a constant field $Y$ by a function $f\in\cc$ produces a constant vector field $fZ$, i.e. $\gg$ is
naturally defined as a Lie algebra over $\cc$.

In $\gg$ we distinguish an ideal $\gg_{||}\subset \gg$ formed by the fields, which 
have the form $Y=f\partial_t$ on the complement
of a compact, where $f\in \cc$.
Such fields will be called {\it parallel}.

The complement to $\gg_{||}$ in $\gg$ is formed by perpendicular vector fields.
A constant vector field $Y$ on $M$ will be called {\it perpendicular}, if the support of the function
$Y(t)$ is compact. The set of all perpendicular vector
fields $Y$ on $M$ will be denoted $\gg_{\perp} \subset \gg$. It is a Lie subalgebra of $\gg$.

Any $Z\in\gg$ can be uniquely represented as
$Z=  Y + f\partial_t, \quad\text{where}\quad Y\in \gg_\perp.$
In this representation $f=Z(t)$, $f\in \cc$.

We will denote by $A(\gg)$ the associative algebra of differential operators generated by $\gg$.
We will call $D\in A(\gg)$ {\it differential operators with constant coefficients}. Note that such operators
are in fact arbitrary on the compact parts of $M$. It is clear that any
$D\in A(\gg)$ can be represented as a finite linear combination with coefficients in $\cc$ of the operators of the form
$D = \partial_t^nY_1Y_2\dots Y_m,\quad\text{where}\quad Y_1,\dots, Y_m\in \gg_\perp.$

\subheading{1.3. Functions of polynomial growth} Let $M$ be a manifold with cylindrical structure.
Suppose that $f:M\to \R$ is a smooth function. We will say that $f$ {\it has polynomial growth} if for any differential
operator with constant coefficients
$D\in A(\gg)$ there exist constants $n>0$ and $C>0$, such that
$$|D(f)(x)| \le C\cdot t(x)^n\tag1-3$$
for all $x$ in a COCS.
The set of all functions of polynomial growth on $M$ will be denoted $\p(M)$. It is clearly a ring.
Any differential operator $Y\in A(\gg)$ maps $\p(M)$ into itself. Functions of the form 
$\sum_{i=0}^m t^i a_i$, where $a_i\in \cc$, are clearly examples of functions of polynomial growth.
As another example we mention the function $a\sin t$, where $a\in \cc$.

We will now define the
space $\po^r(M)$ of polynomial $r$-forms on $M$ as the set of all smooth $r$-forms $\omega$ on $M$ with real values
such that $\omega(X_1, X_2,\dots , X_r)$ belongs to $\p(M)$ for any $r$-tuple of constant vector fields
$X_1, X_2, \dots, X_r\in \gg$.

The exterior derivative $d$ maps $\po^r(M)$ into $\po^{r+1}(M)$.
This follows from the formula
$$
\aligned
&d\omega(X_0, \dots, X_{r}) = \\
&\sum_{j=0}^{r} (-1)^j X_j(\omega(X_0, \dots,\hat X_j,\dots, X_{r}))
+ \\ &\sum_{i<j}(-1)^{i+j}\omega([X_i,X_j], X_0,\dots,\hat X_i,\dots,\hat X_j,\dots, X_r).
\endaligned
$$
If we assume that the fields $X_0,\dots, X_r$ belong to $\gg$, then each term in the RHS of this
formula belongs to $\p(M)$.

We obtain a chain complex $(\po^\bullet(M),d)$, which we will call {\it the complex of polynomial forms}.

It is clear that the complex of polynomial forms is a subalgebra, i.e.
if $\omega_1, \omega_2\in \po^\bullet(M)$ then $\omega_1\wedge \omega_2\in \po^\bullet(M)$.

The cohomology of the complex $(\po^\bullet(M),d)$ of polynomially growing forms coincides with the De Rham cohomology of $M$; this follows from our Theorem 1 below.

\subheading{1.4. Witten deformation} Let $M$ be a manifold with a cylindrical end.
Suppose that $h: M\to \R$ is a fixed smooth function of polynomial growth $h\in \P(M)$.
Consider the complex
$\po^\bullet (M)$ of polynomial forms on $M$ with the deformed differential
$d_h: \po^i(M) \to \po^{i+1}(M),$
where
$$d_h(\omega) = e^{-h}d(e^h\omega) = d\omega + dh \wedge \omega.$$
Our aim is to compute the homology
of this deformed differential.

An equivalent way of thinking about the 
cohomology of the complex of polynomial forms with respect to
$d_h$ is the following. Consider the subcomplex of the De Rham complex $\Omega^\ast(M)$
consisting of the forms $e^h\omega,$
where $\omega$ is a polynomial form; this subcomplex is invariant under the usual exterior derivative $d$. Its cohomology with respect to $d$ is isomorphic to $H^\ast(\po(M),d_h)$.

It is clear that we may change $h$ in the following way without changing the cohomology
$H^\ast(\po(M),d_h)$. Suppose that $p:M\to \R$ is a function of polynomial growth,
such that $p\ge C\cdot t^{-n}$ holds for some $C>0$ and $n>0$. The last condition guarantees that the function $p^{-1}$ is also polynomially growing. Then the function
$$H= h + \log p$$
produces a gauge equivalent differential
$$d_H\omega = p^{-1}d_h(p\omega), \quad \omega \in \po^\ast(M).$$
For example, we may always put $H=h+\log t +const$ without changing the cohomology.

\subheading{1.5. Bounded functions and fields} We will say that a smooth
function $f:M\to \R$ {\it is bounded with all its derivatives} if for any differential operator
with constant coefficients
$D\in A(\gg)$ the function $D(f):M\to \R$ is bounded.

A vector field $Z$ on $M$ will be called {\it weakly bounded} if it can be represented in the form of a finite
linear combination $Z = \sum f_iY_i$, where $Y_i$ are constant vector fields and the functions $f_i$ are bounded.
We will say that $Z$ is {\it strongly bounded}, if in the above representation the functions $f_i$ are bounded with all their derivatives.

\subheading{1.6. Vector field development} We will compute the cohomology of the complex of
polynomial forms with respect to the Witten
deformation $d_h$ assuming existence of a special vector field $Y$, which we call {\it development for $h$},
cf. below.

Let $c>0$ be a fixed sufficiently large number.
Denote
$$U_c^+ = \{(x,t)\in \T ; h(x,t)>c\},\quad\text{and}\quad U_c^- = \{(x,t)\in \T; h(x,t)<-c\}.\tag1-4$$
We will call $$U_c = {U_c}^+ \cup  {U_c}^-$$
 {\it "the union of remote fibers".}

{\it A development for $h$} is a smooth vector field $Y$ on $M$,
such that for some large $c>0$ the following conditions hold:
\roster
\item"(i)" $Y$ is {\it strongly bounded} (cf. 1.5);
\item"(ii)" for any constant vector field $Z$ the field $t\cdot [Z,Y]$ is weakly bounded;
\item"(iii)" $Y(t)=1$ holds outside of a compact subset;
\item"(iv)" $Y(h/t)\le 0$ on $U_c^-$;
\item"(v)"  $Y(h)\ge 0$ on $U_c^+$.
\endroster

We may think of $Y$ as defining together with the function
$t$ an equivalent (in some sense)
cylindrical structure on $M$, which is somehow adjusted to $h$ (via properties
(iv) and (v)).

Condition (iv) is equivalent to
\roster
\item"(iv')" $Y(h)\le h/t\quad\text{on}\quad U_c^-.$
\endroster

Note that conditions (iv) and (v) are not symmetric. In fact, (iv) is stronger than the requirement
$Y(h)\le 0$ on $U_c^-$ which is analogous to (v).

We see that the function $h$ is decreasing along
the integral curves of $Y$  on the negative remote fibers $U_c^-$ and it is non-decreasing on the
positive remote fibers $U_c^+$. Therefore all remote negative fibers are diffeomorphic.

\subheading{Example } Suppose that $h$ on the tube part $\T=\partial M_0\times \R_+$
(cf. 1.1) has the form $h(x,t) = a(x)t^n$ with $n\ge 1$; here $x\in \partial M_0$ and $t\in \R_+$.
Then the field
$Y=\partial_t$ satisfies all the above requirements.

Now we will state the main results of the paper.

\proclaim{Theorem 1} Let $M$ be a manifold with a cylindrical structure, and let $h:M\to \R$
be a polynomially growing function, admitting a development, cf. 1.6.
Then the cohomology of the complex
of polynomial forms $\po^\bullet (M)$ with respect to the deformed differential $d_h$
is given by
$$H^i(\po^\bullet (M), d_h)\simeq  H^i(M, h^{-1}(\{-c\});\R)\tag1-5$$
for sufficiently large $c>0$. In other words, it is isomorphic to the usual relative
cohomology of the pair ($M$, a remote negative fiber). 
\endproclaim

As a simple application of Theorem 1 we mention the following.

\proclaim{Theorem 2} Let $N$ be a smooth compact manifold and let $g:N\to \R$ be a smooth
function. Consider the product $M=N\times \R$ with the standard cylindrical structure, which is 
defined as follows. Let $(x,\tau)$ denote the coordinates on $M$, where $x\in N$ and $\tau\in \R$; then
$t(x,\tau)= |\tau|$ and $\partial_t=\sign \tau\cdot \partial_\tau$. Define $h:M\to \R$ by $h(x,t)=tg(x)$. Then the cohomology of the complex of polynomially growing
forms on $M$ with respect to the deformed differential $d_h$ is given by
$$H^i(\po^\bullet (M), d_h)\simeq  H^i(N, N-Z)$$
where $Z= g^{-1}(\{0\})$ is the set of zeros of $g$.
In particular, $H^i(\po^\bullet (M), d_h)$ vanishes in all dimensions if $g:N\to\R$ has no zeros.
\endproclaim
\demo{Proof} We consider $N\times [-1,1]=M_0$ as the compact part of $M$. The development
$Y=\partial_t$ exists since $h$ is homogeneous in $t$. The set $U_c^-$
(the union of negative remote fibers) consists of points $(x,t)$ with $t\ge -c/{g(x)}$
(assuming that $g(x)$ is positive) and $t\le -c/{g(x)}$ (assuming that $g(x)$ is negative), cf. Fig 2.
\topinsert
\figsecond
\botcaption{Figure 2}\endcaption
\endinsert
Thus we see that the pair $(M,U_c^-)$ is homotopy equivalent to $(N,N-Z)$.
\qed
\enddemo
Here is another application of Theorem 1.
\proclaim{Theorem 3} Suppose that $M$ is the total space of a vector bundle $\pi : M\to N$
over a compact manifold $N$. Let $h: M\to \R$ be a smooth function, which is a nondegenerate quadratic form on each fiber $\pi^{-1}(p)$ for each $p\in N$. Consider the standard cylindrical
structure on $M$ (cf. example 2 in 1.1) and the complex of polynomially growing forms on $M$.
Then we have
$$H^i(\po (M), d_h)\simeq  H^{i-r}(N, \xi),$$
where $r$ is a nonnegative integer and
$\xi$ is a flat real line bundle over $N$, which are constructed as follows.
We may split the bundle $M$ as the Whitney  sum $M\simeq M_+\oplus M_-$ such that $h$ is
positively (negatively) definite on $M_\pm$; then $r$ is the rank of $M_-$ and $\xi$ is the orientation
bundle of $M_-$.
\endproclaim
\demo{Proof} The union of negative remote fibers in each fiber appears as shown in Fig. 3
\topinsert
\figfourth
\botcaption{Figure 3}\endcaption
\endinsert
Thus we see that the pair $(M,U_c^-)$ is homotopy equivalent to $(D(M_-),\partial D(M_-))$,
where $D(M_-)$ is the unit disk bundle of the negative bundle $M_-$. Applying the Thom isomorphism, cf. \cite{BT},
completes the proof.
\qed
\enddemo

\subheading{1.7. Polynomial forms on $\R^n$}
{\it The standard cylindrical structure on} $\R^n$
is given by the function $r=r(x)= |x|$ and by the radial vector field
$\partial_r=r^{-1}\sum_{i=1}^n x_i\partial_i$ considered on the complement of a ball with center at the
origin. The constant vector fields are linear combinations of $\partial_r$ and of the fields
$X_{ij}= x_j\partial_i-x_i\partial_j,\quad i<j,$
which are perpendicular to the radial lines. Translating our general definitions to the present
situation, we easily see that {\it a functions $f:\R^n\to \R$ is polynomially growing iff
for any sequence $i_1,i_2,\dots, i_k\in \{1,2, \dots, n\}$ there exists $N$ so that}
$$|\partial_{i_1}\partial_{i_2}\dots \partial_{i_k}f(x)|\le |x|^N,\quad \text{for}\quad |x|\ge 1.$$
In other words, $f\in \Cal P(\R^n)$ if and only if all partial derivatives of $f$ admit polynomial
estimates in $r=|x|$. Clearly, $\R[x_1,...,x_n]\subset\Cal P(\R^n)$. The above statement follows since we can represent each field
$\partial_k$ as a linear combination $\partial_k =\sum_{i<j} f_{ij}^k X_{ij} + g^k\partial_r$
of constant vector fields with coefficients $f_{ij}^k, g^k$ which are bounded with all their partial
derivatives.

Similarly, {\it a polynomially growing differential form} $\omega\in \po^k(\R^n)$ is given by
$$\omega = \sum a_{i_1i_2\dots i_k}dx_{i_1}\wedge dx_{i_2}\dots \wedge dx_{i_k},$$
where the coefficients $a_{i_1i_2\dots i_k}$ belong to $\Cal P(\R^n)$.

We want to mention that a similar statement concerning functions having bounded
higher derivatives is only partially correct. Namely, any smooth
function $f:\R^n\to \R$ with the property that $D(f):\R^n\to\R$ is bounded for any differential operator $D$ with constant coefficients as defined in section 1.2 (with respect to the standard cylindrical structure) has bounded partial derivatives $\partial_{i_1}\partial_{i_2}\dots \partial_{i_k}f(x)$. But, the converse does not hold; namely,
there exist functions such that all their partial derivatives are bounded but they are not bounded
with all their derivatives according to our general definition (cf. 1.5). A simple example provides
the function
$f(x,y)=\frac{x}{1+x^2}.$
All its partial derivatives are bounded but the derivative with respect to the constant vector field
$X=y\partial_x-x\partial_y$ equals $y\cdot\frac{1-x^2}{1+x^2}$ and is unbounded.

Here is our main result concerning the polynomials on $\R^n$.

\proclaim{Theorem 4} Let $h\in\R[x_1,...,x_n]$ be a polynomial satisfying either 

A. $h$ is quasi-homogeneous with integral weights $w_1>0, w_2>0,\dots, w_n>0$, having degree 
$\deg_w(h) \ge \max\{w_i\}$, or

B. the Newton polyhedron of $h$ lies under a hyperplane
$\sigma$, and the polynomial $h^{\sigma}$
has no critical points on $\R^n$ except the origin.

Then the cohomology of the complex of 
polynomially growing forms on $\R^n$ with respect to the deformed differential $d_h$
coincides with the cohomology of the remote negative fiber, i.e.
$$H^k(\Cal P\Omega^\bullet (\R^n),d_h)\simeq H^k(\R^n,h^{-1}(-c))
\simeq \tilde H^{k-1}(h^{-1}(-c))\tag1-6$$
for any $k$ and for sufficiently large $c>0$. 
\endproclaim

Recall that the Newton polyhedron of $h$ 
is defined as the convex hull of the points $(i_1,\dots, i_n)$
of the integral lattice 
such that the corresponding monomial appears in $h$ with a nontrivial
coefficient. The notation $h^\sigma$ 
stands for the sum of the terms in $h$, which correspond to the
points lying on $\sigma$.

Note that a generic polynomial $h\in \R[x_1,\dots,x_n]$ of fixed degree satisfies B of Theorem 4.

The proof of Theorem 4 is set forth in section 4. It is based on the following theorem, which we prove in section \S 5.

\proclaim{Theorem 5} Let $w_1>0, \dots, w_n>0$ be the integral weights of the variables $x_1, \dots, x_n$ and let $f:\R^n\to \R$ be a real quasi-homogeneous polynomial  
having no critical points in $\R^n -\{0\}$. Let $g$ be a real polynomial of quasi-degree 
$$\deg_w g \le \deg_w f -\max\{w_i\}\tag1-7$$
and let $h = f+g.$
Then there exists a diffeomorphism $\phi:\R^n \to \R^n$ with the following properties:
\roster
\item"(i)" for any $x\in \R^n -K$, where $K\subset \R^n$ is a compact subset, holds 
$$h(x) = f(\phi(x));$$
\item"(ii)" $\phi$ and $\phi^{-1}$ are given by functions with all partial derivatives 
bounded.
\endroster
\endproclaim

More precisely, in (ii) we mean the following. Suppose that $\phi$ is 
given by $\phi(x_1, \dots, x_n) = (\phi_1(x_1, \dots, x_n), \dots, \phi_n(x_1, \dots, x_n))$. Then 
all the functions
$$\frac{\partial^{\alpha_1+\dots+\alpha_n} \phi_i}{\partial x_1^{\alpha_1}\dots
\partial x_n^{\alpha_n}}: \R^n \to \R$$
are bounded. Similarly with respect to $\phi^{-1}$.

Note that in the case $w_1= \dots =w_n=1$ the condition (1-7) is equivalent to $\deg g < \deg f$.

\heading{\bf \S 2. Equivalence relations between cylindrical structures}\endheading

The main question we address in this section is the following. 
Let $M$ denote a fixed manifold with cylindrical structure $(t, \partial_t)$. 
Suppose we produce a
different cylindrical structure with the same function $t$ and with another field $\partial_t$;
we are interested to know when the two cylindrical structures
\roster
(A)  have the same class of functions bounded with all derivatives,

(B)  have the same class of polynomially growing functions.
\endroster
It is clear that (A) implies (B). Correspondingly, we have two different equivalence relations between
the cylindrical structures with given function $t$. Our main results in this section (Theorems 2.5
and 2.7) give necessary conditions for (A) and (B).

We may describe these results as follows. Let $(t,\partial_t)$ be the given cylindrical structure,
and let $Y_\perp$ be a perpendicular vector field on $M$. Then we obtain the deformed cylindrical
structure $(t,\partial_t +Y_\perp)$. If one only assumes that $Y_\perp$ is strongly bounded then
(as easy examples show) the functions, which are bounded with all derivatives with respect to 
$(t,\partial_t )$, may have exponentially growing higher derivatives with respect to the other
structure $(t,\partial_t +Y_\perp)$. 
In Theorem 2.5 we assume that the field $Y_\perp$ decays as $t^{-1}$ and we show
that (B) holds. In Theorem 2.7 we have a stronger assumption that the field 
$Y_\perp$ decays as $t^{-1-\epsilon}$ (where $\epsilon >0$) and we show
then that (A) holds.

\proclaim{2.1. Lemma} A smooth vector field $Z$ is strongly bounded if and only if all
the Poisson brackets of the form
$[X_1,[X_2,[\dots [X_k,Z]]]\dots ]$
are weakly bounded for any set of constant vector
fields $X_1,\dots, X_k$ .
\endproclaim

\demo{Proof} In one direction the Lemma is trivial. Namely, if $Z$ is strongly bounded then all the
brackets with constant vector fields are weakly bounded.
To show the converse we observe that we may find a system of finitely many smooth functions
$\phi_i:M\to \R$, where $i=1,\dots, N$, so that
\roster
\item $\sum\phi_i =1$ holds on a COCS;
\item $\partial_t(\phi_i)=0$ on a COCS; in other words  $\phi_i\in \cc(M)$;
\item the support of $\phi_i$ is sufficiently "thin", which means that for any $i=1, 2,\dots, N$
there are finitely many constant vector fields $Y_1, \dots, Y_n$ (where $n=\dim M$),
commuting with each other $[Y_i, Y_k]=0$, so that
for any point $p\in \supp \phi_i$ with $t(p)$ large enough the vectors
$Y_1(p), \dots, Y_n(p)$ form a basis of the tangent space $T_p(M)$.
\endroster
It follows that for any smooth vector field $Z$ on $M$ we have $\phi_i Z = \sum_{j=1}^n f_j Y_j,$
where $f_j$ are uniquely determined and smooth. Note that the functions $\phi_i$ are bounded
with all their derivatives. Given a field $Z$ as above the following conditions are equivalent:
\roster
\item
$Z$ is weakly (strongly) bounded;
\item all the fields $\phi_i Z$, where $i=1, \dots, N$, are weakly (strongly) bounded;
\item for any index $i=1, \dots, N$, the corresponding functions $f_j$ (where $j=1, \dots, n$)
are bounded (bounded with all their derivatives).
\endroster
It follows that
$$[Y_{1}, [Y_{2},\dots, [Y_{k},\phi_i Z]\dots ] = Y_1(Y_2(\dots Y_{k}(f_j))\dots ) Y_j$$
and thus we obtain (if the assumptions of the Lemma are satisfied) that all the functions
$Y_{1}(Y_{2}(\dots Y_{k}(f_j))$ are bounded; the latter means that the functions $f_j$ are bounded with all their derivatives. \qed
\enddemo

In the next Lemma we will use the following notions. We will say that a smooth vector field $Z$ is
{\it weakly polynomial} if the field $t^{-n}Z$ is weakly bounded for some $n\ge 0$.
We will say that a smooth vector field $Z$ is
{\it strongly polynomial} if it can be represented as a finite linear combination $\sum f_i X_i$, where $X_i$ are
constant vector fields and the functions $f_i$ are polynomially growing, cf. 1.3.

\proclaim{2.2. Lemma} A smooth vector field $Z$ is strongly polynomial if and only if all
the Poisson brackets of the form
$[X_1,[X_2,[\dots [X_k,Z]]]\dots ]$
are weakly polynomial for any set of constant vector
fields $X_1,\dots, X_k$ .
\endproclaim
\demo{Proof} It is similar to the proof of Lemma 2.1 above. \qed
\enddemo

\proclaim{2.3. Lemma} Let $g:M\to M$ be a diffeomorphism obtained as the time-one map of a strongly bounded vector field $Z$ on $M$.
Then for any strongly bounded vector field $Y$ on $M$ the field $g_\ast(Y)$ is strongly bounded.
Also, for any function $f:M\to \R$, which is bounded with all its derivatives, the function $g^\ast(f) = f\circ g$ is bounded with all its derivatives. Similarly, if $f:M\to \R$ is polynomially growing, then $g^\ast(f)$ is polynomially growing.
\endproclaim

Recall that the value of the field $g_\ast(Y)$ at a point $m\in M$ is defined as $dg(Y_{g^{-1}(m)})$.

\demo{Proof} Consider first the case of functions on the Euclidean space $\R^n$. The
coordinates on $\R^n$ will be denoted $x_1, \dots, x_n$.
We also denote $\frac{\partial}{\partial x_j}$ by $\partial_j$.

We will use the following terms. We say that a smooth function
$f:\R^n\to \R$ {\it is bounded with all its derivatives} if the function
$D(f)$ is bounded for any operator $D$ of the form
$D=\partial_{i_1}\partial_{i_2}\dots \partial_{i_r}$, where $i_1,\dots, i_r\in\{1,2,\dots,n\}$.
In other words, all partial derivatives of $f$ are bounded as functions on $\R^n$.

We will also use the notion of {\it a polynomially growing function} $f:\R^n\to \R$. These are functions
such that for any operator $D$ of the form
$D=\partial_{i_1}\partial_{i_2}\dots \partial_{i_r}$, where $i_1,\dots, i_r\in\{1,2,\dots,n\}$,
the function $D(f)$ can be majorized by $|x|^N$ for some $N$, i.e.
$$|D(f)(x)|\le |x|^N,\quad |x|\ge 1.$$

Now, let $Z=\sum_{i=1}^n h_i\partial_i$
be a vector field, such that the functions $h_1, \dots, h_n$ are bounded with all their derivatives.
Let $g: \R^n\to \R^n$ be the diffeomorphism determined as the time-one map of $Z$. We claim
that {\it for any function $f:\R^n\to \R$, which is bounded with all its derivatives, the function
$g^\ast(f)= f\circ g$ is also bounded with all its derivatives}. Similarly, {\it if $f$ is
polynomially growing then $g^\ast(f)$ is polynomially growing}. Since the partial derivative
$\partial_{i_1}\partial_{i_2}\dots \partial_{i_r}(f(g(x))$ can be expressed as the polynomial in
the partial derivatives of $f$ and of $g$ we obtain that all the above statements follow if we will show
that the functions of the form
$$\partial_{i_1}\partial_{i_2}\dots \partial_{i_r}g_j(x),\quad r\ge 1, \quad j\in \{1,2,\dots,n\}
\tag2-1$$
are bounded, where $g(x)= (g_1(x), \dots, g_n(x))$.

Now, let $g^\tau(x),\quad \tau\in \R$, denote the one-parameter group of diffeomorphisms $g^\tau:\R^n\to \R^n$
determined by the field $Z$, so that $g(x)=g^1(x)$ and $g^0(x) =x$. We have
$g^\tau(x) = (g_1^\tau(x),\dots, g^\tau_n(x))$. Since
$$\frac{d}{d\tau}(\partial_{i_1}\partial_{i_2}\dots \partial_{i_r}g_j^\tau (x))=
(\partial_{i_1}\partial_{i_2}\dots \partial_{i_r}h_j)(g^\tau(x))\tag2-2$$
and since we know that the partial derivatives of $h_j$ are bounded, we obtain that the
functions (2-1) are bounded. Thus the Lemma follows in the case of $\R^n$.

\topinsert
\figthird
\botcaption{Figure 4}\endcaption
\endinsert

Consider now the case of an arbitrary manifold $M$ with a cylindrical end.
The cylindrical structure on $M$ will be denoted $(X,t)$. We will reduce the general
case to the case of $\R^n$ considered above by finding an embedding $i:M\to \R^n$ for some large
$n$ so that the composite
$M\buildrel i\over \to \R^n\buildrel {x_n}\over \to \R$
coincides with $t:M\to \R$ outside a compact subset $K\subset M$ (see Fig. 4).
Also we will assume that the vector field $i_\ast(X)$ coincides with $\partial_n=\frac{\partial}{\partial_n}$. We will show that for any vector field $Z$ on $M$
we can construct a vector field $\tilde Z$ on $\R^n$ so that $\tilde Z|_{M}=Z$ and $\tilde Z$ is
bounded (i.e. it can be represented as $\sum_{i=1}^nh_i\partial_i$ so that the functions $h_i:\R^n\to \R$ are bounded with all their derivatives) if $Z$ is bounded. Also we will show that any function
$f:M \to \R$ can be canonically extended to a smooth function $\tilde f:\R^n\to \R$ and the extension
$\tilde f$ is bounded, or polynomially growing (in the sense explained in the beginning
of the proof) if and only if the initial function $f:M\to \R$ is bounded, or polynomially growing, correspondingly. This will clearly prove the Lemma.

We fix a small tubular neighborhood $U\supset M$ and let $\pi: U\to M$ be the projection. We assume
that $U$ is the union of straight line intervals of length $2\epsilon$, which are orthogonal to $M$
and intersect $M$ at the center.
Let $\chi:U \to \R$ be a smooth function, which is identically 0 near $\partial U$ and identically
1 near $M\subset U$. We will also assume that $\chi(x_1,\dots,x_n)$ is independent of $x_n$ for
$x_n\ge c$; in other words, the partial derivative $\partial_n(\chi)$ has a compact support.

The following provides an extension procedure for functions.
For any smooth function $f:M\to \R$ we set
$$\tilde f (x) = \chi(x)f(\pi(x)),\tag2-3$$
where $x\in \R^n$. This defines a smooth function on $\R^n$, having support in $U$.

Suppose now that $Z$ is a vector field on $M$. We may write $Z=\sum_{i=1}^nh_i\partial_i$,
where $h_i:M \to \R$. We extend these functions $h_i$ as explained above  to get $\tilde h_i:\R^n\to \R$ and set
$$\tilde Z= \sum_{i=1}^n \tilde h_i\partial_i.\tag2-4$$
This gives an extension procedure for vector fields.

{\it Given $i_1,\dots, i_r\in \{1,\dots, n\}$ there exists a differential operator
$D$ so that for any smooth function $f:M\to \R$,
$$\partial_{i_1}\partial_{i_2}\dots \partial_{i_r}\tilde f = D(f)\tag2-5$$
and $D$ has the form
$$D(f) =\sum h_{j_1\dots j_k}\cdot (Y_{j_1}\dots Y_{j_k}(f)\circ \pi),\tag2-6$$
where $Y_1, \dots, Y_m$ is a fixed finite set of constant vector fields on $M$, the indices $j_1, j_2, \dots, j_k$ run
through the set $\{1,2,\dots, m\}$, the sum is finite, $k\le r$,
and the functions $h_{j_1\dots j_k}:\R^n\to \R$ are smooth having support in $U$ and are independent
of the last variable $x_n$ for large $x_n\ge c$.}

This follows by induction from a similar statement for the special case $r=1$, which claims that
$$\partial_i\tilde f = h_\emptyset\cdot (f\circ \pi )+ \sum_{j=1}^mh_j\cdot (Y_j(f)\circ \pi),\tag2-7$$
where the functions $h_\emptyset, h_1, \dots, h_m$ are smooth, having support in $U$,
and are independent of the last variable
$x_n$ for large $x_n\ge c$. To show this one observes that there exist finitely many
constant vector fields $Y_1, \dots, Y_m$ on $M$ and smooth functions
$\mu_j:U\to \R$, $j=1, \dots, m$, which are independent
of the last variable $x_n$ for large $x_n\ge c$,
so that
$$(d\pi)_p(\partial_i) = \sum_{j=1}^m \mu_j(p)\cdot Y_j|_{\pi(p)},\qquad p\in U.\tag2-7$$
Then (2-7) holds with $h_\emptyset =\partial_i\chi$ and $h_j=\chi\cdot\mu_j$.

Now we may show that {\it if $f:M\to \R$ is bounded with all its derivatives
(or polynomially growing) then $\tilde f:\R^n\to \R$ is bounded with its all derivatives
(or polynomially growing, correspondingly).} Indeed, assuming that $f$ is bounded
with all the derivatives, we see using formulae (2-5) and (2-6) that
any partial derivative of $\tilde f$ is bounded; here we use the fact that the smooth
functions $h_{j_1\dots j_k}$ have support in $U$ and are independent of the last variable $x_n$
for $x_n\ge c$.

Now we want to prove the inverse statement,  claiming that $f$ is bounded with all its derivatives (or is polynomially
growing) if the function $\tilde f$ has the corresponding property.
With this aim we observe:

{\it If $Z$ is a constant vector field on $M$, then $\tilde Z$ is bounded.}
Indeed, if $Z$ equals $X=\partial_t$, then $\tilde Z= \chi\partial_n$ outside a compact subset,
which is bounded.
If $Z$ is a constant perpendicular vector field
then on $M$ we have $Z=\sum_{i=1}^{n-1}h_i\partial_i$, where the functions $h_i$ do not depend on the
last coordinate $x_n$ for large $x_n\ge c$; therefore the functions $\tilde h_i:\R^n\to \R$
are also independent of $x_n$ for $x_n\ge c$.
But they have a compact support with respect to the other variables.

Assume now that $f:M\to \R$ is such that $\tilde f: \R^n\to \R$ is bounded with all its derivatives.
If $Z_1,\dots, Z_r$ are some constant vector fields on $M$ then the fields $\tilde Z_j$ are bounded
and the function $\tilde Z_1\dots \tilde Z_r(\tilde f)$ is bounded. Hence
$$\tilde Z_1\dots \tilde Z_r(\tilde f)|_M= Z_1\dots  Z_r(f)\tag2-9$$
is also bounded, which proves that $f$ is bounded with all its derivatives.

The analogous statement for polynomial growth follows similarly.

This completes the proof. \qed

\enddemo

\proclaim{2.4. Lemma} Let $M$ be a manifold with a cylindrical structure $(t, \partial_t)$ and
let $Y$ be a vector field on $M$ having the form,
$$Y=Y_\perp +\partial_t,\tag2-10$$
where the field $Y_\perp$ is perpendicular (i.e. $Y_\perp(t)=0$ holds on a COCS).
Assume that $Y_\perp$ is strongly bounded and for
any constant vector field $X$ the field $t[X,Y_\perp]$ is weakly bounded.
Then
\roster
\item"(a)" any smooth vector field on $M$, satisfying $[Z,Y]=0$ on a COCS is strongly polynomial,
i.e. it can be represented as a finite linear combination $\sum f_iX_i$, where $X_i$
are constant vector fields and $f_i$ are polynomially growing.
\item"(b)" Any constant vector field $X$ on $M$ can be represented as a finite linear
combination $\sum f_i Z_i$, where $Z_i$ are fields satisfying $[Z_i,Y]=0$ on a COCS
and $f_i$ are polynomially growing.
\endroster
\endproclaim

\demo{Proof} We will use a similar method of embedding into $\R^n$ as in the proof of Lemma 3.
Namely, fix an embedding $i:M\to \R^n$, so that
$M\buildrel i\over \to \R^n\buildrel {x_n}\over \to\R$
coincides with $t:M\to \R$ outside a compact subset $K\subset M$.
Also we will assume that the vector field $i_\ast(X)$ coincides with $\partial_n=\frac{\partial}{\partial_n}$. We will denote by $U$ a small tubular neighborhood of
$M$, which is obtained as the union of straight intervals of length $2\epsilon$ orthogonal to $M$.
As above, we will denote by $\pi: U\to M$ the projection. $\chi:\R^n\to \R$ will denote a smooth
function which is identically 0 near $\partial U$ and identically 1 near $M$.

Given the field $Y=Y_\perp +\partial_t$ on $M$, we may use
the extension construction described in the proof of Lemma 3,
to construct the vector field $\tilde Y_\perp$ on $\R^n$.
On the domain $x_n\ge c$ it has the form
$$\tilde Y_\perp=\sum_{j=1}^{n-1} h_j\partial_j,\tag2-11$$
where $c>0$ is some fixed number.
 From the proof of Lemma 2.3 we know that $\tilde Y_\perp$ is bounded and so all the derivatives
$\partial_{i_1} \partial_{i_2}\dots \partial_{i_r}(h_j)$ are bounded;  also all the functions
of the form
$$x_n\cdot \partial_{i}(h_j)\tag2-12$$
are bounded.

Let $\hat\partial_n$ denote a bounded vector field on $\R^n$, which coincides with
$X=\partial_t$ on $M$, and with $\partial_n=\frac{\partial}{\partial x_n}$ on the domain $x_n\ge c$.
We set
$$\hat Y = \tilde Y_\perp + \hat\partial_n.\tag2-13$$

We claim now that
{\it any vector field $\hat Z$ on $\R^n$, commuting with $\hat Y$ on the
domain $x_n\ge c$ and having a compact support for $x_n\le c$, is polynomially growing.}
Indeed, if
$g^\tau:\R^n\to \R^n$, $\tau \in\R$,
denotes the flow determined by $\hat Y$, then, since $[\hat Z,\hat Y]=0$ for $x_n\ge c$,
we obtain that
$$\hat Z_p = (dg^\tau)_\ast (\hat Z_{g^{-\tau}(p)}),\qquad p\in \R^n,\quad x_n(p)\ge c,\tag2-14$$
where $\tau = x_n(p)-c.$
If $g^\tau (p) = (g^\tau_1(p),\dots, g^\tau_n(p))$, and
$\hat Z =\sum \hat Z_i\partial_i$, then
$g^\tau_n(p) = x_n(p) +\tau$
and from (2-14) we obtain
$$\hat Z_i(p) = \sum_{j=1}^n (\frac{\partial g^\tau_i}{\partial x_j}\cdot \hat Z_j)|_{g^{c-x_n(p)}(p)}.\tag2-15$$
Also, for $i=n$ we get
$$\hat Z_n(p) = \hat Z_n(g^{c-x_n(p)}(p)),\quad x_n(p)\ge c.\tag2-16$$
We want to show that {\it all functions of the form
$$\partial_{i_1}\partial_{i_2}\dots \partial_{i_r} (\hat Z_j)\tag2-17$$
admit polynomial estimates in the variable $x_n$}. Note that in the RHS of
formulae (2-15) and (2-16) only the restriction of $\hat Z_j$ onto
the set $x_n=c$ are used, where these functions have compact support.
Therefore our statement would follow from (2-15) and (2-16) if we establish that
{\it given a sequence $i_1, i_2,\dots, i_r\in \{1,2,\dots, n\}$, there exist $c'>0$ and $l>0$
so that for all $\tau\le 0$ and $p\in \R^n$,}
$$|\partial_{i_1}\partial_{i_2}\dots \partial_{i_r}(g^\tau_j(p))|\le c'\cdot (\max\{x_n (p)-c,0\})^l.
\tag2-18$$
In fact, we have the differential equations
$$\frac{d}{d\tau}(g^\tau_j(p)) = h_j(g^\tau_1(p), \dots, g^\tau_n(p)),
\quad j=1,2,\dots, n-1\tag2-19$$
and we view the last variable $g_n^\tau(p)=x_n(p)+\tau$ as a parameter. By differentiating we obtain
$$ \frac{d}{d\tau}(\partial_i(g^\tau_j(p))) = \sum_{k=1}^n
\frac{\partial h_j}{\partial x_k}|_{g^\tau(p)} \cdot
\partial_i g^\tau_k|_{p}.\tag2-20$$
Applying Lemma 2.6 below we obtain, using our assumption
$$|\frac{\partial h_j}{\partial x_k}(x_1, x_2,\dots, x_n)|\le \frac{c''}{x_n},\tag2-21$$
that the first derivatives $\partial_i(g^\tau_j(p))$ admit polynomial
estimates of the form (2-18).
We will now use induction. We have
$$ \frac{d}{d\tau}(\partial_{i_1}\dots \partial_{i_r}(g^\tau_j(p))) = \sum_{k=1}^n \frac{\partial h_j}{\partial x_k}|_{g^\tau(p)}  \cdot
\partial_{i_1}\dots \partial_{i_r} g^\tau_k(p) + \Cal L,\tag2-22$$
where $\Cal L$ denotes the terms involving the
lower order derivatives, about which we already know the polynomial estimates
by induction. Applying Lemma 2.6 again gives our statement.

Now we may complete the proof of statement (a) of the Lemma.
If $Z$ is a smooth vector field on $M$ so that
$[Z,Y]=0$ on a COCS, then we may clearly construct an extending field $\hat Z$ on $\R^n$,
so that $[\hat Z, \hat Y]=0$ holds for $x_n\ge c$ and $\hat Z|_M =Z$. We may achieve this by,
first, extending the given field $Z$ somehow onto the hyperplane $x_n=c$ and, secondly,
by translating this field using the flow $g^\tau$ determined by $\hat Y$ which gives a field
commuting with $\hat Y$. From the above description it is clear that we will have that
$\supp(\hat Z) \cap \{x_n=c\}$ is compact.

If $X_1, X_2, \dots, X_k$ are some constant vector fields on $M$ then
$\tilde X_1, \tilde X_2, \dots, \tilde X_k$ (constructed as in the proof of Lemma 2.3)
are strongly bounded fields on $\R^n$ and therefore
$[\tilde X_1, [\tilde X_2, \dots ,[\tilde X_k,\hat Z]]]\dots ]$
has a polynomial estimate in $x_n$. Hence we obtain that for some $m>0$ the field
$t^{-m}[X_1, [X_2, \dots ,[X_k,\hat Z]]]\dots ]$ is weakly bounded. Applying Lemma 2.2
completes the proof of (a).

Now we want to prove (b). For any $i=1, 2,\dots, n$ we will denote by $\hat \partial_i$
the unique vector field on $\R^n$, satisfying $[\hat \partial_i, Y]=0$ for $x_n\ge c$ and
$\hat \partial_i|_{x_n=c}=\chi \partial_i$. We may write
$$\hat \partial_i=\sum_{j=1}^n s_{ij}\partial_j,\tag2-23$$
where $s_{ij}$ are uniquely determined smooth functions with support in $U$; we
have shown before that they are polynomially growing. We claim now that for their determinant
$s = \det(s_{ij})$
an estimate of the form
$$|s(x_1,\dots, x_n)|\ge a\cdot |x_n|^{-l}\tag2-24$$
holds for some constants $a>0, l>0$. Together with the  above established fact,
that $s_{ij}$ and all their derivatives
have polynomial estimates in $x_n$, this would clearly imply that we can solve system (2-23) to obtain
$$\partial_i = \sum k_{ij}\hat \partial_{j}\tag2-25$$
and {\it the coefficients $k_{ij}$ and all their partial derivatives will admit polynomial estimates
in the variable $x_n$.}

To prove (2-24) we observe that the function $s$ satisfies the differential equation
$$\frac{ds}{dx_n} = \alpha\cdot s,\quad\text{where}\quad \alpha=\sum \partial_i h_i.\tag2-26$$
Since we know from the assumptions of the Lemma
that $|\alpha|\le \frac{b}{x_n}$ with some $b>0$ we obtain
$$
\aligned
&|s(x_1,\dots,x_{n-1},x_n)| = |s(x_1,\dots,x_{n-1},c)| \cdot \exp[\int_c^{x_n} \alpha(x_1,\dots,
x_{n-1},\xi)d\xi]\ge \\
&\ge |s(x_1,\dots,x_{n-1},c)|\cdot [x_n/c]^{-b}.
\endaligned
\tag2-27$$
This proves (2-24) since the intersection of 
the support of $s$ with the hyperplane $x_n=c$ is
compact.

To finish the proof of (b), assume that $X$ is 
a constant vector field on $M$. Then we may
write $\tilde X= \sum_{i=1}^n \alpha_i \partial_i$, 
where $\alpha_i$ does not depend on the last variable 
$x_n$ for $x_n\ge c$. Hence we have
$$\tilde X= \sum_{j=1}^nf_j\hat \partial_j,\tag2-28$$
where
$f_j=\sum_i\alpha_i k_{ij}$ and all the partial 
derivatives of $f_j$ have polynomial estimates on
$x_n$. One may construct the following family of 
$n\times n$-matrices depending on $x\in M$:
\roster
\item $A(x)$ smoothly depends on $x$;
\item $A(x)$ is constant along the trajectories 
of the field $Y$, i.e. $Y(a_{ij})=0$ where
$A(x)=(a_{ij}(x))$;
\item $A(x)$ is an idempotent $A(x)^2=A(x)$, 
and for any $x\in M$ the image $\{A(x)\cdot v;v\in \R^n\}$ equals
$T_x(M)$, the tangent space to $M$.
\endroster
We may choose such a family arbitrarily for $x_n=c$ and then extend it using (2).
Having such $A(x)$ we will define the following fields on $M$
$Z_j = A(x)\cdot \hat\partial_j$.
Thus $[Z_j,Y]=0$ on a COCS and $X=\sum f_jZ_j$.
This completes the proof. \qed

\enddemo

\proclaim{2.5. Theorem} Let $M$ be a manifold with a cylindrical structure $(t, \partial_t)$ and
let $Y$ be a vector field on $M$ having the form,
$$Y=Y_\perp +\partial_t,\tag2-29$$
where the field $Y_\perp$ is perpendicular (i.e. $Y_\perp(t)=0$ holds on a COCS).
Assume that $Y_\perp$ is strongly bounded and for
any constant vector field $X$ the field $t[X,Y_\perp]$ is weakly bounded.
We may view the pair $(Y,t)$ as defining another cylindrical structure on $M$. Then
the spaces of polynomial differential forms $\po^\bullet(M)$ with respect to both structures $(t,\partial_t)$ and $(t,Y)$ coincide.
\endproclaim
\demo{Proof} This clearly follows from Lemma 2.4.\qed
\enddemo

\proclaim{2.6. Lemma} Let $x: [1,\infty) \to \R^N$ be a 
smooth curve satisfying the differential equation
$$\frac{dx}{d\nu} = A(\nu)\cdot x(\nu) + b(\nu),$$
where $A(\nu)$ is an $N\times N$-matrix and $ b(\nu)$ 
is an $N$-vector, depending continuously
on the parameter $\nu\in [1,\infty)$ and satisfying
$$|A(\nu)|\le \frac{k}{\nu}, \qquad |b(\nu)|\le l\cdot\nu^k$$
for some constants $k, l  >0$. Assume also that $|x(1)|\le m$.
Then 
$$|x(\nu)|\le \nu^{k}(l\cdot\nu+m-l)\tag2-30$$
for all $\nu\ge 1.$
\endproclaim
\demo{Proof} Let $r=r(\nu)$ denote $|x(\nu)|$. Then we have
$$\frac{dr}{d\nu} \le |\frac{dx}{d\nu}|\le |A(\nu)|\cdot r+|b(\nu)|\le \frac{k}{\nu}\cdot r + l\nu^k. \tag2-31$$
Thus we obtain
$$\frac{d}{d\nu}(\frac{r}{\nu^k})\le l,$$
$$\frac{d}{d\nu}(\frac{r}{\nu^k}-{l}\cdot\nu )\le 0$$
$$\frac{r}{\nu^k}-l\cdot\nu\le m-l,$$
$$r\le \nu^{k}(l\cdot\nu+m-l)$$
and the Lemma follows. \qed
\enddemo

\proclaim{2.7. Theorem}Let $M$ be a manifold with a cylindrical structure $(t, \partial_t)$ and
let $Y$ be a vector field on $M$ having the form,
$$Y=Y_\perp +\partial_t,\tag2-32$$
where the field $Y_\perp$ is perpendicular.
Assume that for some $\epsilon >0$ the field
$t^{1+\epsilon}Y_\perp$
is strongly bounded.
Then we have two cylindrical structures $(\partial_t, t)$ and $(Y,t)$ and we claim that they are equivalent
in the following sense: any smooth function $f:M\to \R$ is bounded with all its derivatives  with respect to one of these
cylindrical structures if and only if it is bounded with all its derivatives with respect to the other.
\endproclaim
\demo{Proof} The proof is similar to the proofs of Lemma 2.4 and Theorem 2.5. Instead of Lemma 2.6 we will use
Lemma 2.8 (cf. below). It follows by induction applied to the system of equations (2-22),
where the term $\Cal L$
involves  (nonlinearly) some lower order partial derivatives of $g^\tau_k(p)$ (which are bounded  by induction)
and the higher order derivatives of $h_j$, which enter linearly and
have estimates by $k\cdot |x_n|^{-1-\epsilon}$.
\qed
\enddemo

\proclaim{2.8. Lemma} Let $x: [1,\infty) \to \R^N$ be a smooth curve satisfying the differential equation
$$\frac{dx}{d\nu} = A(\nu)\cdot x(\nu) + b(\nu),$$
where $A(\nu)$ is an $N\times N$-matrix and $ b(\nu)$ is an $N$-vector, depending continuously
on the parameter $\nu\in [1,\infty)$ and such that
$$|A(\nu)|\le \frac{k}{\nu^{1+\epsilon}}, \qquad |b(\nu)|\le \frac{k}{\nu^{1+\epsilon}}$$
for some constants $\epsilon >0$, $k>0$ and all $\nu\ge 1$. Assume also that 
$|x(1)|\le m$.
Then $x(\nu)$ is bounded and, moreover, 
$$|x(\nu)|\le (m+1)e^{k/\epsilon}\tag2-33$$
for all $\nu\ge 1.$
\endproclaim
\demo{Proof}As in the proof of Lemma 2.6 we obtain for $r=|x(\nu)|$ the inequalities
$$r'\le \frac{k}{\nu^{1+\epsilon}}\cdot r + \frac{k}{\nu^{1+\epsilon}},$$
$$\frac{d}{d\nu}[(r+1)e^{k/\epsilon{\nu^\epsilon}}] \le 0$$
and thus
$$(r+1)e^{k/\epsilon {\nu^\epsilon}}\le (m+1)e^{k/\epsilon}$$
and the result follows.\qed
\enddemo

\proclaim{2.9. Lemma} Let $M$ be a manifold with cylindrical structure and let $Y$ be a development for $h:M\to \R$. Suppose that $\phi: M\to \R$ is a smooth function
such that
\roster
\item $e^{-h}Y(e^h\phi)=\psi: M\to \R$ is
of polynomial growth;
\item for some $c>0$ the function $\phi$ and all its partial
derivatives with respect to the constant vector fields on $M$ are polynomially growing on $U_c^-$;
this means that for any constant vector fields $X_1, X_2, \dots, X_k$ on $M$ there are constants
$C>0$ and $m>0$ such that
$|X_1X_2\dots X_k(\phi)(x,t)|\le C\cdot (t^m +1)$
for all $(x,t)\in U_c^-$.
\endroster
Then $\phi$ is polynomially growing.
\endproclaim

\demo{Proof} The development $Y$ together with the function $t$ define another
cylindrical structure on $M$ in the sense described in Theorem 2.5.
We think of $M$ as represented as the union $M=M_0\cup \T$ of the compact part $M_0$ and the tube
part $\T$ (as in 1.1) and the tube part $\T$ is represented as the product
$\partial M_0\times \R_+$, so that the trajectories of $Y$ are the curves $t\mapsto (x,t)$.
Here $x\in \Gamma = \partial M_0$.
We will consider functions $f:M\to \R$ restricted to the tube part as functions of $x\in \Gamma$ and $t$.

We will prove Lemma by induction. First, we assume that
$\psi$ satisfies an inequality of the form
$$|\psi|\le Ct^n +C',\tag2-34$$
and show then that
$\phi$ satisfies a similar inequality
$|\phi|\le C_1t^{n+1} +C_2.$

Solving the equation $Y(e^h\phi) = e^h\psi$, we have
$$\phi(x,t) = e^{h(x,t_0)-h(x,t)}\phi(x,t_0)+\int_{t_0}^t e^{h(x,\tau)-h(x,t)}\psi(x,\tau)d\tau.\tag2-35$$
Here $t_0$ is any fixed number, for example $t_0=1$.
Assuming that $(x,\tau)\in U_c^-$ for sufficiently large $\tau$ we have instead
$$\phi(x,t) = -\int_{t}^{\infty} e^{h(x,\tau)-h(x,t)}\psi(x,\tau)d\tau.\tag2-36$$

We will represent $\Gamma$ as the union $\Gamma_+\cup\Gamma_-\cup \Gamma_0$, where
$\Gamma_\pm$ is the set of all points $x\in \Gamma$, such that $(x,\tau)\in U_c^\pm$ for large
$\tau$, and $\Gamma_0$ is the set of all points such that $|h(x,\tau)|\le c$ for all $\tau$. Note that the flow determined by $Y = \partial_t$ preserves the sets $U_c^+$ and $U_c^-$.
For $x\in \Gamma_\pm$ we will denote by $T(x)$ the smallest value of $t$ such that $(x,t)$ belongs to
$U_c^\pm$.

Assuming that $x\in \Gamma_0$ or $x\in \Gamma_\pm$ and $t<T(x)$ we have from (2-35)
$$|\phi(x,t)|\le e^{2c}(|\phi(x,t_0)|+\int_{t_0}^t |\psi(x,\tau)|d\tau).\tag2-37$$

If $x\in \Gamma_+$ and $t\ge T(x)$ then we obtain using (2-35)
$$|\phi(x,t)|\le |\phi(x,t_0)|+\int_{t_0}^t |\psi(x,\tau)|d\tau.\tag2-38$$

Now we want to find a similar estimate for $x\in\Gamma_-$ and $t\ge T(x)$.
We will need the following identity
$$
\int_t^\infty e^{-a\tau}\tau^nd\tau = e^{-at}\cdot t^n\cdot n!\cdot a^{-1}\cdot \sum_{i=0}^n \frac{1}{i!}(ta)^{i-n},\qquad a>0.\tag2-39
$$
If $x\in \Gamma_-$ and $t\ge T(x)$, then for $\tau \ge t$ we have
$$h(x,\tau) -h(x,t)\le \frac{h(x,t)}{t}(\tau -t) \le -\frac{c}{t}\cdot (\tau -t).\tag2-40$$
Hence
$$|\phi(x,t)|\le \int_t^\infty e^{h(x,\tau)-h(x,t)}|\psi(x,\tau)|d\tau \le
 \int_t^\infty e^{-c\cdot t^{-1}\cdot (\tau -t)}|\psi(x,\tau)|d\tau.\tag2-41$$
Therefore, using (2-39) with $a=c/t$, we see that
in all the above cases (2-37), (2-38), (2-41),  if $\psi$ satisfies an inequality of the form 
$|\psi|\le Ct^n +C'$, then
$\phi$ satisfies an inequality of the form
$|\phi|\le C_1t^{n+1} +C_2.$

Suppose now that $Z$ is a vector field commuting with $Y$. Then from $Y(e^h\phi) = e^h\psi$ we obtain
$Y(e^hZ(\phi)) = e^h\tilde \psi,$
where
$\tilde \psi = Z(\psi) - YZ(h)\cdot \phi.$
Therefore, if $Z$ is any polynomial vector field commuting with $Y$, we obtain that
$$|Z(\phi)|\le C\cdot t^m +C',$$
assuming that $\psi$ and $Z(\psi)$ satisfy similar inequalities.

The Lemma now follows by induction. We use here Theorem 2.5 stating that any constant vector field
$Z$ on $M$ (with respect to the original cylindrical structure)
can be represented as a finite sum $\sum f_iZ_i$, where $Z_i$ are polynomial vector fields commuting with
$Y$ (i.e. constant fields with respect to the new cylindrical structure)
and $f_i$ are functions of polynomial growth.\qed
\enddemo

\heading{\bf \S 3.  Proof of Theorem 1}\endheading

Because of the obvious isomorphisms
$$H^i(M,U_c^-;\R) \simeq H^i(M, h^{-1}(-c);\R),$$
Theorem 1 would follow if we establish an isomorphism
$$\alpha: H^j(\po^\bullet (M), d_h) \to H^j(M, U_c^-). \tag3-1$$
We understand the relative homology $H^j(M,U_c^-)$ as the homology of the cone
of the chain map
$r: \Omega^\bullet(M) \to \Omega^\bullet(U_c^-).$
Here $\Omega^\bullet(M)$ and $\Omega^\bullet(U_c^-)$ denote the De Rham complex formed by smooth forms on $M$
and $U_c^-$. The chain map $r$ is the restriction of forms.
Recall that by the definition $\Cone^j(r) = \Omega^j(M)\oplus
\Omega^{j-1}(U_c^-)$, and the boundary homomorphism of the cone acts as follows:
$$d(\omega, \omega') = (d_M\omega, -d_U\omega' +r(\omega))\tag3-2$$
for $\omega\in \Omega^j(M)$ and $\omega'\in\Omega^{j-1}(U_c^-).$

Suppose that $\omega\in \po^j(M)$. We assume that $M$ is represented as in (1-1).
On the tube part $\Cal T$ we may write
$$\omega = \omega_{||}\wedge dt + \omega_\perp,\tag3-3$$
where $\omega_{||}$ and $\omega_\perp$ are $(j-1)$-form and $j$-form on $\T$,
which can be viewed as forms on $\Gamma=\partial M_0$, depending on $t$.
We will call them {\it the parallel and the perpendicular components} of $\omega$.
Clearly $\omega_{||}$ and $\omega_\perp$ are polynomial iff $\omega$ is.

Define the following $(j-1)$-form $\omega'$ on $U_c^-$:
$$\omega' (x,t) = (-1)^j\cdot \int_t^\infty e^{h(x,\tau)}\omega_{||}(x,\tau)d\tau.\tag3-4$$
We view the form $e^h\omega_{||}$ as a form on $\Gamma$ depending on the parameter $\tau$, and the integral
above is understood as integrating a curve in the space $\Omega^{j-1}(\Gamma)$. The integral
converges since $\omega$ is assumed to be polynomial and $h(x,\tau) \le -a\tau$ on $U_c^-$ for some $a>0$
(depending on $x$).

We claim now that the map $(\po^\bullet (M),d_h) \to \Cone(r),$
given by
$$\omega \mapsto (e^h\omega, \omega')\tag3-5$$
is a chain map; here $\omega'$ is given by (3-4).
The {\it proof} reduces to check of the following equality
$$-d\omega' +e^h\omega = (-1)^{j+1}\int_t^\infty (e^hd_h\omega)_{||}d\tau \quad\text{over}\quad U_c^-,\tag3-6$$
which follows easily from the identity
$$
(e^hd_h\omega)_{||} = (de^h\omega)_{||} =
d_\Gamma(e^h\omega_{||}) +
(-1)^j\frac{d(e^h\omega_\perp)}{dt}.\tag3-7
$$

We will denote by
$\alpha: H^\ast(\po^\bullet (M),d_h)\to H^\ast(M,U_c^-)$ the induced map on the cohomology

Next we prove that $\alpha$ {\it is a monomorphism.}
Suppose that $\omega\in \po^j(M)$ is a polynomial form with $d_h\omega =0$, such that its image under $\alpha$ bounds.
This means that there exist a $(j-1)$-form $\nu\in \Omega^{j-1}(M)$ and a $(j-2)$-form
$\nu'\in \Omega^{j-2}(U_c^-)$ such that
$$d\nu = e^h\omega\quad\text{and}\quad \omega' =\nu - d\nu'\quad\text{over} \quad U_c^-.\tag3-8$$
Let $c'$ be any number less than $c$. Choose a smooth function $\chi: M\to \R$ such that $\chi|_{U_{c'}}=1$ and
$\chi|_{M-U_c^-}=0$. Then denoting $\mu = \nu -d(\chi\nu)$ we have
$$d\mu = e^h\omega\quad\text{and}\quad \mu|_{U_{c'}}= \omega'= (-1)^j\int_t^\infty e^h\omega_{||}d\tau.\tag3-9$$

Now we want to show that we may change $\mu$ so that it behaves in a special way on the
tube part $\T$. Namely, write
$$\mu = \mu_{||} \wedge dt + \mu_\perp\tag3-10$$
on $\T$. Then we have
$$
d_\Gamma(\mu_{||}) + (-1)^{j-1}\frac{d\mu_\perp}{dt} = e^h\omega_{||}, \quad\text{and}\quad
d_\Gamma(\mu_\perp) = e^h\omega_\perp.\tag3-11
$$

We claim that {\it we may choose $\mu$ so that additionally to (3-9) we will have
$\mu_{||}(x,t) =0$ for large $t$}. To do so we choose a smooth function $\phi:M\to \R$ so that
it is identically zero on $M_0$ and identically 1 on $\Gamma\times[1+\epsilon,\infty)$. Then we set
$$\mu' = \mu - d((-1)^{j+1}\phi\cdot\int_1^t\mu_{||}(x,\tau)d\tau).\tag3-12$$
Note that $d\mu' = e^h\omega$ and $\mu'|_{U_{c'}}= \mu|_{U_{c'}}= \omega'$ since ${\mu_{||}|}_{U_{c'}}=0$
as is seen from (3-9). Moreover, from (3-12) we obtain that $\mu'_{||}(x,t)=0$ for $t>1+\epsilon$.

Now we have
$$\frac{d\mu'_\perp}{dt} = (-1)^{j-1}e^h\omega_{||},\quad \mu'_\perp|_{U_{c'}^-}=
(-1)^j\int_t^\infty e^h\omega_{||}d\tau.\tag3-13$$
We want to apply Lemma 2.9 to conclude that the form $e^{-h}\mu'$ is polynomially growing.
This would show that the initial form $\omega$ is cohomologous to zero in
$(\po^\bullet (M),d_h)$, since $\omega= d_h(e^{-h}\mu')$.
In fact in order to apply Lemma 2.9 we have to know that the form
$$\int_t^\infty e^{h(x,\tau)-h(x,t)}\omega_{||}(x,\tau)d\tau\tag3-14$$
is polynomial on $U_c^-$. It follows using inequalities (2-40) and identity (2-39) similarly to the
proof of (2-41).

We want to show that $\alpha: H^\ast(\po^\bullet (M),d_h)\to H^\ast(M,U_c^-)$ is an epimorphism.
In fact we will show below that {\it any cohomology class $\xi \in H^j(M,U_c^-)$ may be
represented by a smooth closed $j$-form
$\omega$ so that $\omega|_{U_c^-}=0$ and it behaves on the tube part as follows:
$\omega_{||}=0$ and $\omega_\perp$ does not depend on $t$ for large $t$. }
Then it follows that $\alpha$ is onto since the form $e^{-h}\omega$ belongs
to $\po^j(M)$ and its cohomology class with respect to $d_h$ is mapped onto $\xi$ by $\alpha$.

It is clear that we may realize $\xi$ by a closed form $\omega$ on $M$ with $d\omega=0$
and $\omega|_{U_c^-}=0$.
Then on the tube part
$\omega = \omega_{||}\wedge dt + \omega_\perp$
and we have
$$d_\Gamma(\omega_{||})+(-1)^j\frac{d}{dt}(\omega_\perp) =0, \quad
d_\Gamma(\omega_\perp) = 0.\tag3-15$$
Let $\phi: M\to \R$ denote a smooth function which is identically 0 on $M_0$ and identically 1
on $\Gamma\times [1+\epsilon,\infty)$.
The form
$$\omega'= \omega - d((-1)^{j+1}\phi\cdot \int_1^t\omega_{||}dt)\tag3-16$$
is cohomologous to $\omega$ and clearly satisfies
$\omega'|_{U_c^-}=0,\quad \omega'_{||}=0$ where
$\omega'_\perp$ does not depend on $t$ for large $t$.
This completes the proof.\qed

\heading{\bf \S 4. Proof of Theorem 4 }\endheading

\subheading{4.1. Case A: quasi-homogeneous polynomials}
Here we prove Theorem 4 in case A. 
We will construct a development, and so the result
follows from Theorem 1.
 
Let $w_1>0, w_2>0,\dots, w_n>0$ be the weights of $x_1, x_2,\dots, x_n$
and let $h=h(x_1, \dots, x_n)$ be a quasi-homogeneous 
polynomial with these weights of degree $d>0$.
In other words for each monomial $x_1^{i_1}\dots x_n^{i_n}$, 
which appears in $h$, 
$\sum w_k i_k =d.$
We assume that
$d\ge w_i\quad \text{for each}\quad i=1, \dots, n$.
Define $\rho =(\sum w_kx_k^2)^{1/2}$ and set
$$Y= \frac{r}{\rho^2}\sum_{j=1}^n w_jx_j\partial_j.$$
Then, as it is easy to see,
$$Y(r) =1,\quad\text{and}\quad Y(h)=\frac{r}{\rho^2}d\cdot h.$$
Since
$$\frac{r^2}{\rho^2}d\ge 1,$$
we obtain that the sign of the function $Y(h)$ coincides with the 
sign of $h$ and the inequality 
$Y(h)\le h/r$ holds for $h\le 0$ (compare 1.6, inequality (iv')). 
Therefore, $Y$ will provide a development for $h$ assuming that the 
field $rY_\perp$, where
$Y_\perp = Y - \partial_r,$
is strongly bounded. 

One has
$$rY_\perp = \sum_{j=1}^n\frac{w_jr^2-\rho^2}{\rho^2}x_j\partial_j=
\sum_{i<j}\frac{x_ix_j(w_i-w_j)}{\rho^2}X_{i,j},$$
where $X_{i,j}$ denotes the constant field $X_{i,j}=x_i\partial_j - x_j
\partial_i$.
Since the coefficients $\frac{x_ix_j(w_i-w_j)}{\rho^2}$ are constant 
along the radial lines,
we obtain that 
the field $rY_\perp$ is actually a constant vector field (cf. 1.2).
This completes the proof in case A. \qed

\subhead 4.2. Case B: polynomials with quasi-homogeneous leading 
form\endsubhead
Suppose that under the conditions of Theorem 4, Case B, the hyperplane $\sigma$ consists of the
lattice points $(i_1, \dots, i_n)$ satisfying the equation
$$w_1i_1 + w_2i_2 +\dots + w_ni_n=d,$$
where $w_i>0$ are integers.
Let $\delta$ denote $\max \{w_i\}$. We will consider the numbers $w_1, \dots, w_n$ as 
the weights of the 
variables $x_1, \dots, x_n$. $h^\sigma$ is the sum of all monomials $q$ in $h$ with $\deg_w (q) =d$.

Suppose first that $h$ contains no monomials $q$ with nonzero coefficients such that
$d-\delta <\deg_w(q) < d$. Then statement B follows from Theorem 5 and case A of Theorem 4.

We may reduce the general statement B to the special case considered in the previous paragraph
as follows. Consider the diffeomorphism
$\Phi: \R^n \to \R^n$ given by $\Phi(y_1, \dots, y_n)=(x_1, \dots, x_n)$, where
$$x_i = y_i^k + y_i,\qquad i=1, 2, \dots, n.$$
Here $k$ is a fixed odd positive integer with $k \ge \delta$. We observe that $\Phi$ and $\Phi^{-1}$
are given by polynomially growing functions. Therefore, we obtain that 
$\Phi^\ast: \po^\bullet(\R^n)\to \po^\bullet(\R^n)$ is an isomorphism and 
$$d_h(\Phi^\ast\omega) =\Phi^\ast(d_{\tilde h}\omega),\quad \omega \in \po^\bullet(\R^n),$$
where $\tilde h = h\circ \Phi$. It is clear that the Newton diagram of $\tilde h$ will contain 
no monomials
$q$ with $kd -k<\deg_w(q)<kd$, and so we may apply to $\tilde h$ 
the special case considered above.
\qed

\heading{\bf \S 5. Proof of Theorem 5}\endheading

Consider the family of polynomials $h_\tau =f+(1-\tau)g$ where $\tau\in [0,1]$.
Then $h_0=h$ and $h_1=f$. We are going to find a continuous family of diffeomorphisms
$\Phi^\tau$ with $\Phi^0 = \id$, so that $h(x) = h_\tau(\Phi^\tau(x))$ for all $x$ outside a compact subset.

We will assume that the integral weights $w_1>0, \dots, w_n>0$ are fixed and denote
$$\<x\> = \sqrt{x_1^{2/w_1}+\dots+x_n^{2/w_n}}.\tag5-1$$

From our assumptions on $f$ it follows that there exist constants $a>0$ and $r_1>0$ so that for all 
$\tau \in [0,1]$ and all $x$ with $\<x\>\ge r_1$,
$$\sum_{j=1}^n \<x\>^{2w_j}(\frac{\partial h_\tau}{\partial x_j})^2 \ge a\cdot \<x\>^{2d},\tag5-2$$
where $d=\deg_w f$. Indeed, the term of the highest $w$-degree in the LHS of (5-2) is
$$\sum_{j=1}^n \<x\>^{2w_j}(\frac{\partial f}{\partial x_j})^2;$$
it is quasi-homogeneous of degree $2d$ and is strictly positive in a neighborhood of infinity.

Consider now the following system of ordinary differential equations on $\R^n$:
$$\frac{dy_i}{d\tau} = v_i(y_1, \dots, y_n, \tau),\quad i=1, 2, \dots, n\tag5-3$$
with the initial conditions
$$y_i(0) = x_i,$$
where the functions $v_i$ are given by
$$v_i(y, \tau) = \rho(y)\frac{g(y)\cdot\<y\>^{2w_i}\cdot\frac{\partial h_\tau}{\partial y_i}(y)}
{\sum_{j=1}^n [\<y\>^{2w_j}\cdot \frac{\partial h_\tau}{\partial y_j}(y)^2]},$$
and $\rho:\R^n\to [0,1]$ is a smooth function, vanishing identically in $\<y\>\le r_1$,
and $\rho(y)$ is identically 1 outside the set $\<y\> \ge r_1+1$. Inequality (5-2)
guarantees that $v_i$ are well defined. 

Also, 
since $g$ has degree $\le d-\delta$, where $\delta$ denotes $\max\{w_i\}$, we obtain that for some $b>0,$
$$|g(y)|\le b\<y\>^{d-\delta}\tag5-4$$
for all $y$ with $\<y\>\ge 1$. Similarly, for some $c>0$,
$$|\frac{\partial h_\tau}{\partial y_i}| \le c\cdot \<y\>^{d-w_i}\tag5-5$$
for all $y$ with $\<y\>\ge 1$ and all $i$. This shows that
$$|v_i(y,\tau)|\le \frac{bc}{a^2}\cdot \<y\>^{w_i-\delta}\le \frac{bc}{a^2}$$
for $\<y\>\ge r_1+1$ (since we assume that $w_i\le \delta$ for all $i$)
and therefore the solution to the system (5-3) exists. 

Now we observe: {\it there exists a constant $R>0$ such that for any solution $y(\tau)$ of (5-3), 
$\tau\in [0,1]$, }
$$\<y(\tau)-y(0)\>\le R.\tag5-6$$

Another crucial property of our system (5-6) is the following:
{\it 
if $y(\tau)$ is a solution to (5-3)  such that $\<y(0)\>\ge r_1+1+R$, then $h_\tau(y(\tau))$ is independent of
$\tau$ and therefore }
$$f(y(1)) = h(y(0)).\tag5-7$$

We may now define 
our diffeomorphism $\phi: \R^n\to \R^n$ as the time one map of the system (5-3).
In other words, 
$\phi(x) = y(1)$
where $y(\tau)$ is a solution to (5-3) with the initial condition $y(0) = x$. (5-7) shows that property (i)
is already satisfied. To finish the proof we have to show that all the partial derivatives of $\phi$ and
 $\phi^{-1}$ are bounded.

Consider first the partial derivatives $\frac{\partial \phi_i^\tau}{\partial x_k}(x)$.
Here $\phi^\tau:\R^n \to \R^n$, where $\phi^\tau(x) = (\phi^\tau_1(x),\dots, \phi^\tau_n(x))$
is the time $\tau$ diffeomorphism determined by the system (5-3), 
$\tau\in[0,1]$, $\phi^0 =\id$, $\phi^1=\phi$. These derivatives satisfy the following system of linear
differential equations
$$\frac{d}{d\tau}(\frac{\partial \phi_i^\tau}{\partial x_k}) = 
\sum_{j=1}^n \frac{\partial v_i}{\partial x_j}(\phi_i^\tau(x))\cdot \frac{\partial \phi_j^\tau}{\partial x_k}\tag5-8$$
with the initial conditions
$$\frac{\partial \phi_i^0}{\partial x_k}=\delta_{ik}.\tag5-9$$
It is easy to see that for $\<x\> >r_1 +1$,
$$|\frac{\partial v_i}{\partial x_j}(x)|\le c_{ij}\cdot \<x\>^{w_i-\delta- w_j},\tag5-10$$
where $c_{ij}$ is a constant. Since $w_i-\delta- w_j < 0$, we obtain that given $\epsilon >0$
we may find $r_2> r_1+1+R$ large enough so that for all
$x\in \R^n$ with $\<x\>\ge r_2$ and all $\tau\in [0,1]$,
$$|\frac{\partial \phi_i^\tau}{\partial x_k}-\delta_{ik}| <\epsilon.\tag5-11$$
In particular, we obtain that the first derivatives $\frac{\partial \phi_i}{\partial x_k}= 
\frac{\partial \phi_i^1}{\partial x_k}$ are bounded.

Suppose that we have already proven that all partial derivatives of
order less than $p$ of $\phi^\tau$ are bounded.
Consider the derivatives of order $p\ge 2$
$$\frac{\partial^{p} \phi_i^\tau}{\partial x_1^{\alpha_1}\dots
\partial x_n^{\alpha_n}},\qquad \alpha_1+\dots +\alpha_n=p\tag5-12$$
with fixed $\alpha_1, \dots, \alpha_n$.  
They satisfy the system of linear differential equations:
$$\frac{d}{d\tau}(\frac{\partial^{p} \phi_i^\tau}{\partial x_1^{\alpha_1}\dots
\partial x_n^{\alpha_n}}) = 
\sum_{j=1}^n \frac{\partial v_i}{\partial x_j}(\phi_i^\tau(x))\cdot \frac{\partial^{p} \phi_i^\tau}{\partial x_1^{\alpha_1}\dots
\partial x_n^{\alpha_n}} + \psi_i(x)\tag5-13
$$
with initial conditions 
$$\frac{\partial^{p} \phi_i^0}{\partial x_1^{\alpha_1}\dots
\partial x_n^{\alpha_n}} =0,\tag5-14$$
where the functions $\psi_i(x)$ are sums of products involving the higher partial derivatives of the functions $v_i$
and the partial derivatives of $\phi_j^\tau$ of order less then $p$. Since it is clear that all partial derivatives of $v_i$ are bounded, we conclude by induction that
the functions $\psi_i$ are bounded.
Using (5-13) we obtain that the partial derivative (5-12) is also bounded. 
Setting $\tau=1$, proves that all partial derivatives of $\phi$ are bounded. 

The fact that all partial derivatives of $\phi^{-1}$ are also bounded follows similarly using the inequalities
(5-11). \qed

\enddocument